\documentclass[a4paper,12pt]{article}
\usepackage{geometry,latexsym,amssymb,amsmath,color,graphicx}
\usepackage[latin5]{inputenc}
\usepackage{enumerate}
\usepackage{enumitem}
\usepackage[T1]{fontenc}
\usepackage{authblk}
\usepackage[all]{xy}
\usepackage{palatino}
\usepackage{indentfirst}
\usepackage{setspace}
\newtheorem{example}{Example}[section]
\newtheorem{Def}[example]{Definition}
\newtheorem{Exam}[example]{Example}
\newtheorem{Prop}[example]{Proposition}
\newtheorem{Theo}[example]{Theorem}
\newtheorem{Lem}[example]{Lemma}
\newtheorem{Rem}[example]{Remark}
\newtheorem{Cor}[example]{Corollary}
\newenvironment{Prf}{{\bf Proof:} } {\hfill $\blacksquare$
	\mbox{}}

\def\NSGGd/G{\mathsf{NSGGd/G}}
\def\NSCM/(A,B,\alpha){\mathsf{NSCM/(A,B,\alpha)}}
\def\SpPXM/(A,B,\alpha){\mathsf{SpPXM/(A,B,\alpha)}}
\def\LXM/(A,B,\alpha){\mathsf{LXMod/(A,B,\alpha)}}
\def\CXM/(A,B,\alpha){\mathsf{CXMod/(A,B,\alpha)}}
\def\GGdA(G){\mathsf{GpGpdAct(G)}}
\def\GGdC/G{\mathsf{GpGpdCov/G}}

\def\epsilon{\varepsilon}

\begin{document}

\title{Semi-homotopy and semi-fundamental groups}

\author[a]{Ayhan ERCİYES\thanks{Correspondence: ayhan.erciyes@hotmail.com}}
\author[b]{Ali AYTEKİN}
\author[b]{Tunçar ŞAHAN}
\affil[a]{Department of Elementary Mathematics Education, Aksaray University, Aksaray, TURKEY}
\affil[b]{Department of Mathematics, Aksaray University, Aksaray, TURKEY}

\date{}

\maketitle
\begin{abstract} In this study we introduce the notions of semi-homotopy of semi-continuous maps and of semi-paths. We also construct a group structure, which will be called semi-fundamental group, using semi-loops and explore some properties of semi-homotopy and semi-fundamental groups.
\end{abstract}

\noindent{\bf Key Words:} Semi-open sets, semi-closed sets, homotopy, fundamental groups
\\ {\bf Classification:} 54C08, 14F35, 55Q05, 57M05

\section{Introduction}

Homotopy theory studies topological objects up to homotopy equivalence. Homotopy equivalence is a weaker relation than topological equivalence, i.e., homotopy classes of spaces are larger than homeomorphism classes. Therefore, homotopy equivalence plays a more important role than homeomorphism. Homotopy theory is a subdomain of topology. Instead of considering the category of topological spaces and continuous maps, one may prefer to consider as morphisms only the continuous maps up to homotopy. On the other hand the concept of homotopy groups is a way to interpret topological problems to algebraic problems which could be solve much easier. For this reason, homotopy groups, especially fundamental groups, are very powerful tools for this purpose. To obtain further insights on applications of homotopy groups, see for example the books of Brown \cite{Brown} and of Rotman \cite{Rotman}.

The concept of semi-open set in topological spaces was introduced in 1963 by Levine \cite{Lev}. He defined a set $A$ to be semi-open in a topological space if and only if $A$ is between an open subset and the closure of that open. Further, Levine investigated a notion of semi-continuity. After the works of Levine on semi-open sets, various mathematician turned their attention to the generalisations of various concepts of topology by considering semi-open sets instead of open sets. New results are obtained in some occasions and in other occasions substantial generalisations are exibited, by replacing open sets with semi-open sets

In 1971, S. Gene Crossley and S. K. Hildebrand \cite{Cro} introduce semi-closed sets, semi-interior, and semi-closure in a manner analogous to the corresponding concepts of closed sets, interior, and closure. Further, a year later, they  defined that a property of topological spaces is a semi-topological property if there is a semi-homeomorphism which preserves that property \cite{Cro2}. Also, they shown that Hausdorff, separable, and connected properties of topological spaces were semi-topological properties.

S.M.N. Maheshawari and R. Prasad \cite{Mahes} used semi-open sets to define and investigate three new separation axiom called Semi-T$_0$, Semi-T$_1$ and Semi-T$_2$.

Recently, P. Bhattacharyya and B.K. Lahiri \cite{Lahiri} generalised the concept of closed sets to semi-generalised closed sets with the help of semi-openness.

In the light of these works, the main purpose of this paper is to introduce the notions of semi-homotopy and semi-fundamental group using the semi-open sets, to obtain different group structures from topological spaces.

\section{Preliminaries}

The notion of semi-open sets in a topological space was introduced by Levine \cite{Lev} as follows.

\begin{Def}\cite{Lev}
	Let $X$ be a topological space and $A\subseteq X$. $A$ is called \textbf{semi-open}
	provided that there exists an open set $U$ such that $U\subseteq A\subseteq \overline{U}$, where $\overline{U}$ denotes the closure of the set $U$ in $X$.
\end{Def}

Here is a concrete example of semi-open sets.

\begin{Exam}
	Let $\tau =\left\{ X,\emptyset,\{a\},\{a,b\}\right\}$ be the topology on the set $X=\left\{ a,b,c,d\right\}$. Therefore we have semi-open subsets of $X$ as follows:
	\[SO(X)=\{ X,\emptyset,\{a\},\{a,c\},\{a,d\},\{a,b\},\{a,b,c\},\{a,c,d\},\{a,b,d\} \}.\]
\end{Exam}


Following proposition is a well known result for semi-open sets. Hence we omit the proof.

\begin{Prop}\cite{Jen,Lev}\label{unionsemi}
	Union of any collection of semi-open sets in a topological space is also semi-open.
\end{Prop}

\begin{Exam}
	Consider the space of the real numbers with the usual topology. It is easy
	to see that intervals of the form $(a,b),(a,b],[a,b)$ and $[a,b]$ and their
	arbitrary unions are semi-open.
\end{Exam}

\begin{Prop}
	Let $X$ be a topological space and $A\subseteq X$. Then $A$ is semi-open if and only if for each point $x$ in $A$ there exist a semi-open subset $B_x$ of $X$ such that $x\in B_x\subseteq A$.
\end{Prop}

\begin{Prf}
	Let $A$ be a semi-open set in $X$. Thus we can choose the set $B_x$ as $A$ for all $x\in A$.
	
	Conversely assume that for each point $x$ in $A$ there exist a semi-open subset $B_x$ of $X$ such that $x\in B_x\subseteq A$. Then \[\bigcup_{x\in A}B_x=A\] and by Proposition \ref{unionsemi} $A$ is a semi-open subset of $X$.
\end{Prf}

The notion of semi-closedness is introduced in \cite{Cro}. Now we will recall the definition of semi-closed sets and some-properties of semi-closed sets from \cite{Cro}.

\begin{Def}\cite{Cro}
	Let $X$ be a topological space and $C\subseteq X$.
	$C$ is called \textbf{semi-closed} if there exists a closed set $K$
	such that $K^{\circ}\subseteq C\subseteq K$ where $K^{\circ}$ is
	the interior of $K$.
\end{Def}

\begin{Exam}
	Let $\tau =\left\{ X,\emptyset,\{a\},\{a,b\}\right\}$ be the topology on the set $X=\left\{ a,b,c,d\right\}$. Therefore we have semi-closed subsets of $X$ as follows:
	\[SC(X) = \{ X,\emptyset,\{b\},\{c\},\{d\},\{c,d\},\{b,c\},\{b,d\},\{b,c,d\} \}.\]
\end{Exam}

%
%

\begin{Prop}\cite{Cro}
	In a topological space the complement of a semi-open set is semi-closed and vice-versa.
\end{Prop}

Now we will recall the definitions of semi-continuities and some properties of them from \cite{Jen}.

\begin{Def}
	Let $X$ and $Y$ be two topological spaces, $f:X\rightarrow Y$ a function and $p$ a point of $X$. Then $f$ is called
	\begin{enumerate}[label=\textbf{(\roman{*})}, leftmargin=1cm]
		\item \textbf{so-1-continuous} at $p$
		provided for each open set $V$ containing $f(p)$ in $Y$, there exists a
		semi-open set $A$ in $X$ that contains $p$ and $f(A)\subseteq V$,
		\item \textbf{so-2-continuous} at $p$
		provided for each semi-open set $B$ containing $f(p)$ in $Y$, there exists a
		semi-open set $A$ in $X$ that contains $p$ and $f(A)\subseteq B$, and
		\item \textbf{so-3-continuous} at $p$
		provided for each semi-open set $B$ containing $f(p)$ in $Y$, there exists
		an open set $U$ in $X$ that contains $p$ and $f(U)\subseteq B$.
	\end{enumerate}
	If $f$ is so-$i$-continuous at every point of $X$ for a fixed $i$ then $f$ is called \textbf{so-$i$-continuous}.
\end{Def}

Relations between so-$i$-continuous functions, constant functions and continuous functions are given with the following figure.

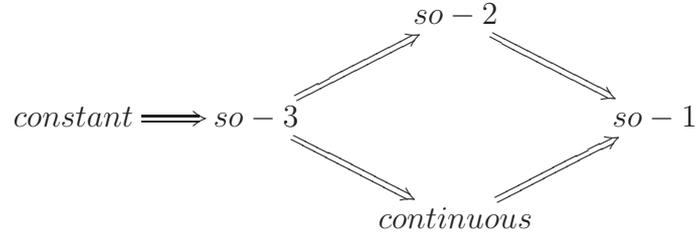
\[\xymatrix{& & so-2 \ar@2{->}[dr] & \\
	constant \ar@2{->}[r] & so-3 \ar@2{->}[ur] \ar@2{->}[dr] & & so-1\\
	& &  continuous \ar@2{->}[ur] &}\]

This figure says that every constant map is so-3-continuous, every so-3-continuous function is both so-2-continuous and continuous, every so-2-continuous function and every continuous function is so-1-continuous.

Following proposition gives a criteria for so-$i$-continuous functions similar to one in classical topology. The proof is also similar, hence we omit.

\begin{Prop}
	Let $X$ and $Y$ be topological spaces and $f:X\rightarrow Y$ a function. Then $f$ is
	\begin{enumerate}[label=\textbf{(\roman{*})}, leftmargin=1cm]
		\item so-1-continuous iff for each open set $V\subseteq Y$, $f^{-1}(V)$ is
		semi-open in $X$,
		\item so-2-continuous iff for each semi-open set $B\subseteq Y$, $f^{-1}(B)$ is
		semi-open in $X$,
		\item so-3-continuous iff for each semi-open set $B\subseteq Y$, $f^{-1}(B)$ is
		open in $X$.
	\end{enumerate}
\end{Prop}

This proposition could be given by using semi-closed sets as follows.

\begin{Prop}
	Let $X$ and $Y$ be topological spaces and $f:X\rightarrow Y$ a function. Then $f$ is
	\begin{enumerate}[label=\textbf{(\roman{*})}, leftmargin=1cm]
		\item so-1-continuous iff for each closed set $K\subseteq Y$, $f^{-1}(K)$ is
		semi-closed in $X$,
		\item so-2-continuous iff for each semi-closed set $M\subseteq Y$, $f^{-1}(M)$ is
		semi-closed in $X$,
		\item so-3-continuous iff for each semi-closed set $M\subseteq Y$, $f^{-1}(M)$ is
		closed in $X$.
	\end{enumerate}
\end{Prop}

so-1-continuous functions are called semi-continuous and so-2-continuous functions
are called irresolute \cite{Cro}. In this paper the unit interval $[0,1]$ will be
denoted by $I$, as a subspace of reel numbers $\mathbb{R}$ with the usual topology.

\begin{Rem}\label{remcomp}
	Let $X$ be a topological space. Then it is easy to see that the identity function $1_X\colon X\rightarrow X$ is so-1-continuous and so-2-continuous but not so-3-continuous. Moreover usual composition of so-2-continuous (resp. so-3-continuous) functions are again so-2-continuous (resp. so-3-continuous). Thus we obtain the category $s$-$\mathsf{Top}$ of topological spaces with morphisms so-2-continuous (irresolute) functions. On the other hand composition of so-1-continuous functions need not to be so-$1$-continuous.
\end{Rem}

\section{Semi-Homotopy}

In this section we will introduce the notions of so-$i$-homotopy of so-$i$-continuous functions, so-2-homotopy type, so-$i$-paths and so-$i$-homotopy of so-$i$-paths, and give some properties. From now on $i$ will symbolize of a fixed element of the set $\{1,2,3\}$ for each item.

\begin{Def}
	Let $X$ and $Y$ be two topological spaces and $f,g\colon X\rightarrow Y$ be two so-$i$-continuous functions. If there exist a function $H\colon X\times I\rightarrow Y$ such that for all $t\in I$ the restrictions of $H$
	\[\begin{array}{rcccl}
	H_t & \colon & X & \longrightarrow & Y \\
	& & x & \longmapsto & H_t(x)=H(x,t)
	\end{array}\]
	are so-$i$-continuous with $H_0=f$ and $H_1=g$, then we say that $f$ and $g$ are \textbf{so-$i$-homotopic}. In this case $H$ is called an \textbf{so-$i$-homotopy} from $f$ to $g$ and this will be denoted by $H\colon f\simeq_i g$ or briefly, by $f\simeq_i g$.
\end{Def}

\begin{Theo}\label{equivrel}
	The relation being so-$i$-homotopic on the set of all so-$i$-continuous functions between two topological spaces is an equivalence relation.
\end{Theo}

\begin{Prf}
	Let $X$ and $Y$ be two topological spaces and $f,g,h:X\rightarrow Y$ be so-$i$-continuous functions.
	
	\begin{enumerate}[leftmargin=2.7cm]
		\item[\textit{\textbf{Reflexivity:}}] If $f:X\rightarrow Y$ define
		\[\begin{array}{rcccl}
		H &\colon & X\times I & \longrightarrow & Y \\
		& & (x,t) & \longmapsto & H(x,t)=f(x)
		\end{array}\]
		for all $x\in X$ and all $t\in I$. It is clear that $F:f\simeq_{i}f$.
		
		\item[\textit{\textbf{Symmetry:}}] Assume that $H:f\simeq_{i}g$, so there is a
		function $H:X\times I\rightarrow Y$ with $H(x,0)=f(x)$ and $H(x,1)=g(x)$
		for all $x\in X$. Define
		
		\[\begin{array}{rcccl}
		G & \colon & X\times I & \longrightarrow & Y \\
		& & (x,t) & \longmapsto & G(x,t)=H(x,1-t)
		\end{array}\]
		
		for all $x\in X$ and all $t\in I$. Since $H$ is so-$i$-continuous,
		\[G_{t}(x)=G(x,t)=H(x,1-t)\] is so-$i$-continuous, and $G_0=g$ and $G_1=f$. Therefore $G:g\simeq_{i}f$.
		
		\item[\textit{\textbf{Transitivity:}}] Assume that $F:f\simeq_{i}g$ and $G:g\simeq_{i}h$. Define
		\[H(x,t) = \begin{cases} F(x,2t), & t\in [0,1/2] \\ G(x,2t-1), & t\in[1/2,1]. \end{cases}\] Therefore $H:f\simeq_{i}h$. Thus $\simeq _{i}$ is an equivalence relation.
	\end{enumerate}
\end{Prf}

Let $X$ and $Y$ be two topological spaces and $f\colon X\rightarrow Y$ be an so-$i$-continuous function. Then the set of all so-$i$-continuous functions from $X$ to $Y$ which are so-$i$-homotopic to $f$ is called the equivalence class (so-$i$-homotopy class) of $f$ and denoted by $[f]_i$.
\[[f]_i=\{g ~|~ g\colon X\rightarrow Y \text{ so-}i\text{-continuous}, f\simeq_i g\}\]

Similar to classical theory, using the new homotopy defined above we will introduce the notion of so-$i$-homotopy equivalence and so-$i$-homotopy type just for the case $i=2$ since the composition of so-2-continuous functions is again so-$i$-continuous.

\begin{Def}
	Let $X$ and $Y$ be two topological spaces. An irresolute function $f\colon X\rightarrow Y$ is called a \textbf{irresolute homotopy  equivalence} if there exist an irresolute function  $g\colon Y\rightarrow X$ such that $gf\simeq_2 1_X$ and $fg\simeq_2 1_Y$. If there is an irresolute homotopy equivalence between two topological spaces then we say that these spaces have the same \textbf{irresolute homotopy type}.
\end{Def}

Now we will give the definition of so-$i$-paths which is the special case of so-$i$-continuous functions. Further we will give a more stronger version of so-$i$-homotopy for so-$i$-paths.

\begin{Def}
	Let $X$ be a topological space, $\alpha\colon I\rightarrow X$ be an so-$i$-continuous function and $\alpha(0)=a$ and $\alpha(1)=b$ . Then $\alpha$ is called an \textbf{so-$i$-path from} $a$ to $b$ in $X$. If $a=b$ then $\alpha$ is called an \textbf{so-$i$-loop} at $a$.
\end{Def}

\begin{Def}\label{semipathcomp}
	Let $\alpha,\beta\colon I\rightarrow X$ be two so-$i$-path in $X$ with $\alpha(1)=\beta(0)$. Then the function
	\[(\alpha\ast\beta) (t) = \begin{cases} \alpha(2t), & t\in [0,1/2] \\ \beta(2t-1), & t\in[1/2,1] \end{cases}\]
	is an so-$i$-path and is called the \textbf{composition} of \textbf{so-$i$-paths} $\alpha$ and $\beta$ in $X$.
	$\alpha\ast\beta$ will be denoted by $\alpha\beta$ for short.
\end{Def}

\begin{Def}
	Let $X$ be a topological space and $\alpha\colon I\rightarrow X$ be an so-$i$-path in $X$. Then the function
	\[\overline{\alpha}\colon I\longrightarrow X\] defined by $\overline{\alpha}(t)=\alpha(1-t)$ is an so-$i$-path
	in $X$ and is called the \textbf{inverse} of $\alpha$.
\end{Def}

\begin{Def}
	Let $X$ be a topological space and $\alpha,\beta\colon I\rightarrow X$ be
	two so-$i$-paths where $\alpha(0)=\beta(0)$ and $\alpha(1)=\beta(1)$. If there
	is an so-$i$-continuous function $F:I\times I\rightarrow X$ such that
	
	\begin{enumerate}[label=\textbf{(\roman{*})}, leftmargin=1cm]
		\item for all $t\in I$ the restrictions of $F$
		\[\begin{array}{rcccl}
		F_t &\colon & I & \longrightarrow & Y \\
		& & s & \longmapsto & F_t(s)=F(s,t)
		\end{array}\]
		are so-$i$-continuous and
		\item $F(s,0)=\alpha (s),$ $F(0,t)=a$, $F(s,1)=\beta (s)$, and $F(1,t)=b$
	\end{enumerate}
	then we say that $F$ is \textbf{so-$i$-homotopy} of \textbf{so-$i$-paths} from $\alpha$ to $\beta$ relative to endpoints
	and denoted by $F\colon \alpha \simeq_i \beta$ rel \^{I}. We will denote this by
	$\alpha \simeq_i \beta$ where no confusion arise.
\end{Def}

\begin{Theo}
	The relation being so-$i$-homotopic relative to endpoints on the set of all so-$i$-paths
	in a topological space is an equivalence relation.
\end{Theo}

\begin{Prf}
	This can be proved by a similar way to the proof of Theorem \ref{equivrel}.
\end{Prf}

\begin{Def}
	Let $X$ be a topological space and $\alpha\colon I\rightarrow X$ an so-$i$-path in $X$.
	Then the set \[[\alpha]_i=\{\beta ~|~ \alpha\simeq_i\beta \ \  \text{ rel \^{I}}\}\] is called
	equivalence class (\textbf{so-$i$-homotopy class}) of $\alpha$.
\end{Def}

\section{Semi-Fundamental groups}

In this section, using the so-$i$-loops, we will construct a group structure on the set of all so-$i$-homotopy classes of so-$i$-loops at a base point of a topological space. Following lemma is a very useful tool to construct this group structure.

\begin{Lem}\label{rho}
	Let $X$ be a topological space, $a,b\in X$ and $\alpha$ be an so-$i$-path from $a$ to $b$. If there is an so-$i$-continuous function $\rho\colon [0,1]\rightarrow[0,1]$ such that $\rho(0)=0$ and $\rho(1)=1$ then $\alpha\rho\simeq_i \alpha$.
\end{Lem}

\begin{Prf}
	First of all note that $\alpha\rho$ is an so-$i$-path from $a$ to $b$. Now we define the so-$i$-homotopy $F\colon \alpha\rho\simeq_i \alpha$ as follows:
	\[\begin{array}{rcccl}
	F &\colon & I\times I & \longrightarrow & X \\
	& & (s,t) & \longmapsto & F(s,t)=\alpha\left((1-t)s+t\rho(s)\right)
	\end{array}\]
	It is easy to see that $F$ is an so-$i$-homotopy from $\alpha\rho$ to $\alpha$.
\end{Prf}

\begin{Prop}\label{propwelldef}
	Let $X$ be a topological space and $\alpha,\beta,\alpha',\beta'\colon I\rightarrow X$ be so-$i$-paths such that $\alpha(0)=\alpha'(0)$, $\alpha(1)=\alpha'(1)=\beta(0)=\beta'(0)$ and $\beta(1)=\beta'(1)$. If $\alpha\simeq_i\alpha'$ and $\beta\simeq_i\beta'$ then $\alpha\beta\simeq_i\alpha'\beta'$.
\end{Prop}

\begin{Prf}
	Let $F$ and $G$ be two so-$i$-homotopy from $\alpha$ to $\alpha'$ and from $\beta$ to $\beta'$, respectively. Then the function $H\colon I\times I\longrightarrow X$ defined by
	\[H(s,t) = \begin{cases} F(2s,t), & s\in [0,1/2] \\ G(2s-1,t), & s\in[1/2,1] \end{cases}\]
	is so-$i$-continuous and defines an so-$i$-homotopy from $\alpha\beta$ to $\alpha'\beta'$.
\end{Prf}

\begin{Prop}\label{propassoc}
	Let $X$ be a topological space and $\alpha,\beta,\gamma\colon I\rightarrow X$ be three so-$i$-paths
	with $\alpha(1)=\beta(0)$ and $\beta(1)=\gamma(0)$. Then \[\alpha(\beta\gamma)\simeq_i(\alpha\beta)\gamma.\]
\end{Prop}

\begin{Prf}
	By the Definition \ref{semipathcomp} compositions $\alpha(\beta\gamma)$ and $(\alpha\beta)\gamma$ are defined as follows:
	\[\alpha(\beta\gamma)(t) = \begin{cases} \alpha(2t), & t\in [0,1/2] \\ \beta(4t-2), & t\in[1/2,3/4] \\ \gamma(4t-3), & t\in[3/4,1] \end{cases}\]
	and
	\[(\alpha\beta)\gamma(t) = \begin{cases} \alpha(4t), & t\in [0,1/4] \\ \beta(4t-1), & t\in[1/4,1/2] \\ \gamma(2t-1), & t\in[1/2,1]. \end{cases}\]
	Now let define a function $\rho\colon I\rightarrow I$ by
	\[\rho(t) = \begin{cases} 2t, & t\in [0,1/4] \\ t+\frac{1}{4}, & t\in[1/4,1/2] \\ \frac{t+1}{2}, & t\in[1/2,1]. \end{cases}\]
	One can see that $\rho$ is an so-$i$-continuous function and $\rho(0)=0$, $\rho(1)=1$. Moreover $(\alpha(\beta\gamma))\rho=(\alpha\beta)\gamma$. Then by Lemma \ref{rho} $\alpha(\beta\gamma)\simeq_i(\alpha\beta)\gamma$.
\end{Prf}

\begin{Prop}\label{propid}
	Let $X$ be a topological space, $x,y\in X$ and $\alpha\colon I\rightarrow X$ be an so-$i$-path from $x$ to $y$. Then
	\[1_x\alpha\simeq_i\alpha\simeq_i\alpha 1_y\] where $1_x$ and $1_y$ are the constant maps at $x$ and $y$, respectively.
\end{Prop}

\begin{Prf}
	First of all let define a function $\rho\colon I\rightarrow I$ by
	\[\rho(t) = \begin{cases} 0, & t\in [0,1/2] \\ 2t-1, & t\in[1/2,1]. \end{cases}\]
	This function satisfies the conditions of Lemma \ref{rho} and $1_x\alpha=\alpha\rho$. Hence $1_x\alpha\simeq_i\alpha$. Similarly by taking $\rho$ as
	\[\rho(t) = \begin{cases} 2t, & t\in [0,1/2] \\ 1, & t\in[1/2,1] \end{cases}\]
	one can show that $\alpha\simeq_i\alpha 1_y$.
\end{Prf}

\begin{Prop}\label{propinv}
	Let $X$ be a topological space, $x,y\in X$ and $\alpha\colon I\rightarrow X$ be an so-$i$-path in $X$ from $x$ to $y$. Then
	\[\alpha\overline{\alpha}\simeq_i 1_x \ \ \text{ and } \ \ \overline{\alpha}\alpha\simeq_i 1_y.\]
\end{Prop}

\begin{Prf}
	Let define a function $F\colon I\times I\rightarrow X$ for all $t\in I$ by
	\[F(s,t) = \begin{cases} \alpha(2s), & s\in [0,t/2] \\ \alpha(s), & s\in[t/2,1-t/2] \\ \alpha(2-2s), & s\in [1-t/2,1]. \end{cases}\]
	This function defines an so-$i$-homotopy from $1_x$ to $\alpha\overline{\alpha}$. Similarly, one can show that $\overline{\alpha}\alpha\simeq_i 1_y$.
\end{Prf}

\begin{Theo}
	Let $X$ be a topological space and $x\in X$. Then the set
	\[\pi_1^i(X,x)=\{[\alpha]_i ~|~ \alpha\colon I\rightarrow X \ \ \text{so-$i$-loop at} \ x\}\]
	of all so-$i$-homotopy classes of so-$i$-loops at $x$ has a group structure with the operation
	\[\begin{array}{rcccl}
	\ast & \colon & \pi_1^i(X,x)\times\pi_1^i(X,x) & \longrightarrow & \pi_1^i(X,x) \\
	& & ([\alpha]_i,[\beta]_i) & \longmapsto & [\alpha]_i\ast[\beta]_i=[\alpha\ast\beta]_i.
	\end{array}\]
\end{Theo}

\begin{Prf}
	Proposition \ref{propwelldef} shows that the operation $\ast$ is well defined. By Proposition \ref{propassoc} the operation is associative. The so-$i$-homotopy class of constant map $1_x$ at $x$ acts as the identity element, i.e. for all $[\alpha]_i\in\pi_1^i(X,x)$
	\[[1_x]_i\ast[\alpha]_i=[\alpha]_i\ast[1_x]_i=[\alpha]_i\]
	by Proposition \ref{propid}. Finally according to Proposition \ref{propinv} for all $[\alpha]_i\in\pi_1^i(X,x)$ the inverse of $[\alpha]_i$ up to the operation $\ast$ is $[\alpha]_i^{-1}=[\overline{\alpha}]_i\in\pi_1^i(X,x)$.
\end{Prf}

This group will be called the \textbf{so-$i$-fundamental group} of $X$ at $x$. In particular $\pi_1^1(X,x)$ will be called \textbf{semi-fundamental group} and  $\pi_1^2(X,x)$ will be called \textbf{irresolute fundamental group}.

\begin{Prop}
	Let $X$ be a topological space, $x,y\in X$ and $\gamma\colon I\rightarrow X$ be an so-$i$-path from $x$ to $y$. Then \[\pi_1^i(X,x)\cong\pi_1^i(X,y).\]
\end{Prop}

\begin{Prf}
	The claimed isomorphism is
	\[\begin{array}{rcccl}
	\gamma_{\star} & \colon & \pi_1^i(X,x) & \longrightarrow & \pi_1^i(X,y) \\
	& & [\alpha]_i & \longmapsto & [\gamma]_i^{-1}\ast[\alpha]_i\ast[\gamma]_i.
	\end{array}\]
\end{Prf}

\begin{Cor}
	In a topological space whose topology is so-$i$-path-connected, i.e. for each pair of elements there exist an so-$i$-path between them, every so-$i$-fundamental group is isomorphic.
\end{Cor}

\begin{Prop}
	Let $s-\mathsf{Top}_{\ast}$ be the category of pointed topological spaces with morphisms so-$2$-continuous (irresolute) functions and $\mathsf{Grp}$ be the category of groups with morphisms group homomorphisms. Then
	\[\begin{array}{rcccl}
	\pi_1^2 &\colon & s-\mathsf{Top}_{\ast} & \longrightarrow & \mathsf{Grp} \\
	& & (X,x) & \longmapsto & \pi_1^2(X,x)
	\end{array}\]
	is a functor.
\end{Prop}

\begin{Cor}
	Let $X$ and $Y$ be two topological spaces. If $f\colon X\rightarrow Y$ is a homeomorphism then
	$\pi_1^2(X,x)\cong\pi_1^2(Y,f(x))$.
\end{Cor}

\section{Conclusion}

It seems that according to these results one can define a more general notion semi-fundamental groupoid following the way in \cite{Brown} and \cite{Rotman}. Further, using the results of the paper \cite{Csaszar} of Cs{\'a}sz{\'a}r it could be possible to develop more generic homotopy types and homotopy groups. Hence parallel results of this paper could be obtained for generalized open sets and for generalized continuity.


\begin{thebibliography}{99}
	\bibitem{Lahiri} Bhattacharyya, P. and Lahiri, B.K., \emph{Semi-generalized closed sets in topology}, lnd. Jr. Math., 29 (1987), 375--382.
	
	\bibitem{Brown} Brown, R., \emph{Topology and groupoids}, BookSurge LLC, North Carolina, 2006.
	
	\bibitem{Csaszar} Cs{\'a}sz{\'a}r, {\'A}., \emph{Generalized open sets}, Acta Mathematica Hungarica, 75(1), (1997), 65--87.
	
	\bibitem{Cro} Crossley, S. and Hildebrand, S.K., \emph{Semi-closure}, Texas J. Sci. 22 (1971), 99--112.
	
	\bibitem{Cro2} Crossley, S. and Hildebrand, S.K., \emph{Semi-topological properties}, Fundamenta Mathematicae 74(3) (1972), 233--254.
	
	\bibitem{Hat} Hatcher, A., \emph{Algebraic topology,} Cambridge University
	Press, 2002.
	
	\bibitem{Lev} Levine, N., \emph{Semi-open sets and semi-continuity in topological spaces}, Amer. Math. Monthly 70 (1963), 36--41.
	
	\bibitem{Mahes} {Maheshawari, S.M.N. and Prasad, R.}, \emph{Some new separation axioms}, Ann. Soco. Sci. Bruxelles 89 (1975), 395--402.
	
	\bibitem{Rotman} {Rotman, J.J.}, \emph{An introduction to algebraic topology}, Springer, 1988.
	
	\bibitem{Jen} Scheers, J.M., \emph{An exploration of semi-open sets in topological spaces}, M.Sc. Thesis, Stephen F. Austin State University, 2011.
	
\end{thebibliography}
\end{document}